 \documentclass[11pt,a4pape,rpdftex]{article}
\usepackage{amsmath, amscd, amssymb, amsthm, latexsym%,cancel
} \usepackage{float}

\usepackage[colorlinks=true,linkcolor=red]{hyperref}
\usepackage{hyperref} 
\newtheorem{theorem}{Theorem}
\newtheorem{proposition}[theorem]{Proposition}
\newtheorem{lemma}[theorem]{Lemma}
\newtheorem{example}[theorem]{Example}

\newtheorem{corollary}[theorem]{Corollary}
\newcommand{\means}[1]{\hbox{$ [\kern -.4em [\, {#1}\, ]\kern -.4em]$}}

\usepackage{showlabels}

\usepackage{showlabels}
\newcommand{\nl}{\medskip \noindent }

\newcommand{\N}{\mathbb{N}}

\newcommand{\gram}{\equiv_{\text{gram}}}

\newcommand{\stel}{\stackrel{*}{\sim}}

\newlength{\transwidth}
\newcommand{\trans}[1]{\settowidth{\transwidth}{$\;_{#1}\;$}
\stackrel{#1}{\overrightarrow{\rule{\transwidth}{0ex}}}}

\begin{document}

\begin{center}
{\large \sc  Grammic monoids with three generators\\
Christian Choffrut,\\
 IRIF, Universit\'{e} Paris Cit\'{e}
\today}\\
\end{center}

\tableofcontents

\bigskip 

{\bf Abstract}Young tableaux are combinatorial objects whose construction can be achieved from words over a finite alphabet  by row or column insertion as shown by Schensted sixty years ago. Recently Abram and Reutenauer studied the action the free monoid on the set of columns by slightly adapting the insertion algorithm. Since the number of columns is finite, this action yields a finite transformation monoid. Here we consider the action on the set of rows. We investigate this infinite monoid in the case of a 3 letter alphabet. In particular we show that it is the quotient of the free monoid relative to a congruence generated by the classical Knuth rules plus a unique extra rule.


 \section{Introduction} 
 
Let $A$ be a finite ordered alphabet.  A  Young tableau
is a finite labeling of the right top quarter discrete plane justified to the left and to the bottom 
in such a way that the rows are nondecreasing from left to right and the columns are increasing from bottom
to top, see Example \ref{ex:young-tableau}. 
Using Schensted's method of inserting a letter in a row (or in a column), a Young tableau can be associated to a word in such a way that all letter occurrences appear once and only once in the tableau. Two word are associated to the   same tableau if and only if they are equivalent in a congruence whose rules were given by the Knuth  \cite[Thm. 6.6]{K70}.
The quotient thus defined is the \emph{plactic monoid} denoted by $\text{Plactic}(A)$.

 Schensted's left insertion of a letter in a column results in a column plus a possible letter.  Ignoring this extra letter yields 
 an action of the free monoid on the (finite) set of columns.
The quotient of
the free monoid by the nuclear congruence is a finite monoid $\text{Styl}(A)$ - the \emph{stylic} monoid, many interesting properties of which are studied in 
 \cite{AR2022}, such as the $J$ classes of these finite monoids of transformation,
the notion of $N$-tableaux as representations of congruence classes in the same 
way Young tableaux represent congruence classes
of the plactic monoid and a minimum set of rules defining the stylic congruence. 
 Here we consider the ``dual'' problem of acting on the set of rows which are finite nondecreasing sequences. The main departure is that this monoid, called \emph{Grammic} monoid as suggested to us by Christophe Reutenauer
   is clearly no longer finite. In  this note we consider a specific aspect, namely the presentation 
 of the monoid when $A$ has three generators. Based on the characterization of the pairs of words defining the same mapping, our main result states that these pairs of words are precisely those that are equivalent in the congruence generated by the Knuth rules and the single new rule $cbab=bcab$ where $a<b<c$. This result is optimal since on $3$ generators the plactic congruence is
 strictly finer than the grammic congruence.
An open issue is to generalize it to an arbitrary number of generators. We conjecture that for $4$ letters the 
 congruence is generated by the Knuth rules and the single new rule $dbac=bdac$ where $a<b\leq c<d$.
 
 I would like to relate the present result with the main Theorem 2.4. of \cite{KO11} investigating 
 the algebra  of the 
 Plactic monoid over a field $K$and which is a continuation of \cite[Theorem 9]{CO02}. It states that the  principal ideal generated by $cbab-bcab$ 
 is one of the two principal ideals  that are minimal  
minimal prime ideals. 

%%%%%%%%%%%%%%%%%%%%%%%%%%%%%%%%%%%%%%%%%%%%%%%%%%%%%%
%%%%%%%%%%%%%%%%%%%%%%%%%%%%%%%%%%%%%%%%%%%%%%%%%%%%%%
\section{Preliminaries} 
%%%%%%%%%%%%%%%%%%%%%%%%%%%%%%%%%%%%%%%%%%%%%%%%%%%%%%
%%%%%%%%%%%%%%%%%%%%%%%%%%%%%%%%%%%%%%%%%%%%%%%%%%%%%%

\subsection{Young tableaux}

The reader is referred to the chapter \cite{LLT02} for an introduction of the plactic monoid
but the basics is recalled now.

Let $A$  be a totally ordered alphabet with $k$ elements
$a_{1}< \cdots < a_{k}$. Given a element $u$ of the free monoid $A^{*}$ 
and a letter $a\in A$ we let $|u|_{a}$ denote the number of occurrences of  $a$ in $u$
and   $|u|$ denote its \emph{length}.
We let the plactic congruence on $A^{*}$  be denoted by $\equiv$
and we recall that it is generated by the Knuth relations \cite[Expression 6.7]{K70}
\begin{equation}
\label{eq:knuth}
\begin{array}{l}
bac\equiv bca \text{ where } a<b\leq c,\\
acb\equiv cab \text{ where } a\leq b< c
\end{array}
\end{equation}

A Young tableau is a labelling by occurrences of letters in $A$
of a lower-order ideal of $\N^{2}$ where the ordering is the natural product ordering on $\N$. 
 Every row is non decreasing from left to right and every column is strictly decreasing from top to bottom.

\begin{example}
\label{ex:young-tableau} We assume $A=\{a,b,c\}$ with $a<b<c$. The following is a Young tableau
 $$
\begin{array}{l}
c   \\
bbc\\
aabb \\
\end{array}
$$

\end{example}

The \emph{height} of a Young tableau is the number of its rows. We will speak of bottom, top rows and 
the like in the natural way. In example \ref{ex:young-tableau} the tableau has height $3$, the bottom row is
$aabb$  and the top row is $c$.

\medskip Schensted gave in \cite{Schen61} an algorithm associating to a word $w$ a Young tableau $P(w)$ 
in which every letter occurrence  appears once and only once 
in such a way that two words $u$ and $v$ are $\equiv$-congruent if they are associated to the same tableau, 
see e.g., \cite[Proposition 6.2.3]{LLT02}. We recall it and reserve henceforth the term ``row'' for a 
nondecreasing sequence of letters and ``column'' for a decreasing sequence of letters.
This algorithm proceeds by insertion of a letter $b$ in a row $b_{1}\cdots b_{p}$. If $b>b_{p}$ 
then the insertion results in the row $b_{1}\cdots b_{p}b$. Otherwise, if $b_{i}$ is the least letter greater than $b$
then the insertion results in the row obtained by substituting $b$ for $b_{i}$. More generally, 
the insertion of $b$ in a Young tableau consists of inserting $b$ in the bottom row. If $b$ is greater 
than the greatest element in the row the procedure stops. Otherwise $b$ bumps the letter 
$b_{i}$  which is inserted in the next upper row 
and the process is iterated as long as necessary.
The Young tableau associated with a word $w$ is obtained incrementally 
by starting from the empty tableau 
and successively inserting the letters of $w$.

Dually the tableau can be constructed by defining how to insert a letter $c$ in a column $c_{1}\cdots c_{p}$.
If $c>c_{1}\cdots c_{p}$ the resulting column is $c c_{1}\cdots c_{p}$. Otherwise if $c_{i}$ is the least letter 
greater than or equal to $c$ (note the difference between row and column insertion) then $c$ is substituted for $c_{i}$
in the column. The construction of the tableau follows the same pattern as in the case of row insertion except that instead of proceeding form bottom to top, it proceeds form let to right.

As a consequence of Schensted construction, there are two specific representatives of an $\equiv$-equivalence class, namely the  \emph{row normal form} 
which is the sequence of the rows of the Young tableau from top to bottom
and the \emph{column normal form} which is the sequence  of the columns of the tableau from left to right.

In example \ref{ex:young-tableau} its row and column normal forms are
$c/bbc/aabb$ and $cba/ba/cb/b$ where we used the backslash symbol
for visual convenience. In particular $cbbcaabb\equiv cbabacbb$.

\subsection{Grammic congruence}

We consider the action of a letter $b$ on a row (we recall that it is a nondecreasing sequence of letters)
which is the insertion of $b$ as explained above in Schensted procedure except that the possible letter which is expelled
is definitely lost. More precisely, we set
$$
a_{1}\cdots a_{n}\cdot b =
\left\{
\begin{array}{ll}
a_{1}\cdots a_{n}\cdot b & \text{if } a_{n}\leq b\\
a_{1}\cdots a_{i-1}ba_{i+1}\cdots a_{n} & \text{if } a_{i} \text{ is the leftmost letter}\\
    & \text{greater than } b,
\end{array}
\right.
$$

The mapping extends to $A^{*}$ and factors through the plactic monoid 
because of the following equalities which can be readily verified.

$$
\begin{array}{l}
a_{1}\cdots a_{n}\cdot  bac= a_{1}\cdots a_{n}\cdot  bca \text{ where } a<b\leq c,\\
a_{1}\cdots a_{n}\cdot  acb = a_{1}\cdots a_{n}\cdot  cab \text{ where } a\leq b< c
\end{array}
$$
We let $u\equiv_{\text{gram}} v$ denote the \emph{grammic} congruence between two words $u$ and $v$ defining the same mappings on the rows and call  \emph{grammic monoid} the quotient of $A^{*}$ by  $\gram$.
The congruence $\gram$ is clearly coarser than the congruence $\equiv$. 
 Alternatively, by identifying $a^{n_{1}}_{1}\ldots a^{n_{k}}_{k}$  with the $k$-tuple 
$(n_{1}, \ldots, n_{k})$ we may consider that the free monoid acts on $\overbrace{\N\times \cdots \times \N}^{ k \text{ times}}$
 in an obvious way.

\begin{lemma}
\label{le:action} Consider the alphabet  $a_{1}, a_{2}, \ldots, a_{k}$ with $k$ letters.
If the mappings are equal then the words have the same commutative image.
\end{lemma}
\begin{proof}

Consider a word $u=\alpha_{1}\cdots \alpha_{n}\in A^{*}$ of length $n$
and the vector $\lambda=(n, n, \ldots, n)\in \N^{k}$. 
By the definition of the action of letter on vectors, for all proper prefixes $w$ of $u$, 
no component of $\lambda\cdot u$ vanishes. Thus

$$
\begin{array}{l}
(\lambda\cdot u)_{i}=
\left\{
\begin{array}{ll}
n+|u_{a_{1}}| & \text{ if } i=1\\
n+|u_{a_{i}}|-|u_{a_{i-1}}|& \text{ if } 1<i\leq k
 \end{array}
 \right.
 \end{array}
$$
Therefore  $(\lambda\cdot u)_{i}=(\lambda\cdot v)_{i}=$
implies $|u|_{a_{i}}= |v|_{a_{i}}$.
\end{proof}

\begin{lemma}
\label{le:bottom-row}
The bottom row of the Young tableau $u$ is $(0,0\ldots, 0)\cdot u$.
\end{lemma}

\begin{proof} Indeed,  we make two preliminary remarks. 
From the very definition of the action of the plactic monoids on the rows, it is clear that 
the image of $(0,0\ldots, 0)$ by the product of 
$p$ columns with the same rightmost letter equal to $a_{i}$, 
is the vector with all components equal to $0$ except component $i$ which is equal to $p$.
Furthermore, the action of $a_{i}$ on a vector in $\N^{k}$  affects 
no components smaller than to $i$. In other words 
for  $u\in (A\setminus \{a_{1}, \ldots, a_{i-1}\})^{*}$ and $(n_{1}, \ldots, n_{k})\in \N^{k}$
there exists  $(m_{i}, \ldots, m_{k})$ such that $(n_{1}, \ldots, n_{k})\cdot u= (n_{1}, \ldots, n_{i-1}, m_{i}, \ldots, m_{k})$.
The claim of the lemma results from the previous observations and the fact that an element of the plactic monoid has a representative as a product
of non decreasing columns  in the following partial ordering over decreasing sequences of letters

\begin{equation}
a_{1}\cdots a_{p}\geq b_{1}\cdots b_{q} \text{ if } p\geq q \text{ and } a_{p-q+1}\geq b_{1},\cdots, a_{p}\geq b_{q}
\end{equation}

(e.g., with Example \ref{ex:young-tableau}  we have 
$cba\geq ba\geq  cb\geq b$). 
\end{proof}

Lemmas \ref{le:action} and \ref{le:bottom-row}  have simple consequences for the case of two- and three-letter alphabets.

\begin{corollary}
\label{cor:two-letters}
Over a two letter alphabet we have
$$
u \equiv v \Leftrightarrow u \equiv_{\text{gram}}  v 
$$
\end{corollary}

\begin{corollary}
I $u \equiv_{\text{gram}}v$ and $u \not\equiv_{\text{plactic}}v$  then  either $u$ or $v$ or both have height at least 3.
\end{corollary}

\begin{proof}
The bottom rows are equal. If they have only  two rows, then the top rows are also equal
\end{proof}

\begin{corollary}
\label{cor:a+c}
Over a three-letter alphabet, if $u \equiv_{\text{gram}}v$ and $u \not\equiv_{\text{plactic}}v$  then  
the row representatives of $u$ and $v$ in the plactic monoid are of the form
\begin{equation}
\label{eq:general-form}
3^{a}2^{b}3^{c}1^{d}2^{e}3^{f} \text{ and } 3^{a'}2^{b'}3^{c'}1^{d'}2^{e'}3^{f'}
\end{equation}
with $a\not=a'$, $0<b=b'$, $e=e'$, $0<d=d'$, $f=f'$, $a,a'\leq b \leq d$,  $b+c, b+c'\leq d+e$ and $a+c=a'+c'$.

\end{corollary}

\begin{proof} By the previous corollary the two words contain an occurrence of each of the three letters and the bottom rows of their Young tableaux are equal, hence $e=e'$, $d=d'>0$ and  $f=f'$. The rows above the bottom one contain the same number of occurrences of $2$. If this number is null then the row normal forms are 
$u=3^{c'}1^{d}2^{e}3^{f} \text{ and } v=3^{c'}1^{d}2^{e}3^{f}$
which implies $c=c'$ because the words are commutatively equivalent and hence $u \equiv_{\text{plactic}}v$.
 Equality $a+c=a'+c'$
is yet another consequence of the commutative equivalence. The remaining inequalities result from
the fact that the expressions of \ref{eq:general-form} are row normal forms.
\end{proof}

\begin{lemma}
\label{le:non-vanishing}
Let $0<i\leq k$ and   $u\in A^{*}$  in row normal form,
$u_{1}/u_{2}/\cdots /u_{p}$.  
Set  $u_{j}=b_{j}w_{j}$ where  $b_{1}>b_{2}>\cdots > b_{p}$ 
are the initial letters of the rows. Consider
$$
m_{i}=\sum_{i<j} \sum_{a<b_{i}} |u_{j}|_{a}
$$
Then for all vectors $x\in \N^{k}$  with $x_{i}> m_{i}$
and for all proper prefixes $w$ of $u$
it holds $(x\cdot w)_{i}>0$.

 \end{lemma}

\begin{proof}

Let $r$ be the greatest integer less than $i$ or $0$ if such an integer does not exist. Since it is clear that $(x\cdot u_{1}u_{2}\cdots u_{p})_{i}=(x\cdot u_{r+1}u_{2}\cdots u_{p})_{i}$ we may assume $r=0$ and compute
$(x\cdot  u_{1}\cdots u_{p})_{i}$.

Rewrite $m_{i}=\sum_{i<j} \mu_{ j} $ and assume $x_{i}> \sum_{i<j} \mu_{j} $. 
Then by induction on $j=1, \ldots, p$ for 
all prefixes 
$v$ of $u_{j}$  we have

$$
(x\cdot u_{1}\cdots u_{j -1}v)_{i}  - \sum_{\ell>j} \mu_{\ell}  >0
$$

\end{proof}

With example \ref{ex:young-tableau} the previous computation 
yields $x_{1}>0, x_{2}>2$ (actually it can be checked that $x_{3}>6$ and thus $(x_{1}, x_{2}+2, x_{3}+6).w=
(x_{1}+|w|_{a}, x_{2}+2 +|w|_{b}-|w|_{a}, x_{3}+6) +|w|_{c}-|w|_{b}$).
\medskip 

In the next lemma given a formal vector 
$X=(\omega_{1}, \ldots,\omega_{i-1}, x_{i}, \omega_{i+1}, \ldots, \omega_{k})$ 
where $\omega_{j}$, $j\not=i$  is an integer constant and $x_{i}$ a variable 
and a word $u\in A^{*}$, we view the expression 
$X\cdot u$ as a function of $\N$ into $\N^{k}$.

\begin{lemma}
\label{le:one-variable}
Let $u,v\in A^{*}$ be two words of length $n$.
If for all $|x_{i}|\leq n+1$ the functions $X\cdot u$ and $X\cdot v$ are equal, then 
these functions are equal for all values of $x_{i}$.
\end{lemma}

\begin{proof}We set  
$
\Delta_{u}(x_{i})= X \cdot u -
X
$ and we show that, viewed as a function of $x_{i}$, 
it is constant for all $x_{i} >|u| $.
 We set $u=b_{1}\cdots b_{n}$ and observe that by Lemma \ref{le:non-vanishing}
 for all prefixes $w=b_{1}\cdots b_{r} $,   $r=0, \ldots, n$,
we have  $(X\cdot w)_{i}\geq n-r$. 
We show by induction that 
 the function 
$\Delta_{b_{1}\cdots b_{r} } (X)$  
 is a constant for all $x_{i}\geq n$. If $r=0$ then $\Delta_{1}(x_{i})$ is the zero vector
 thus we assume $r>0$. We set $w=vb_{r}$ and $b_{r}=a_{j}$ 
 and we compute the action of $v$ on $X$.

\nl  \fbox{  Case 1} $j<i$ and $(X\cdot v)_{h}=0$ for all $j<h<i$ then 
$$
\begin{array}{l}
(X\cdot w)_{j}= (X\cdot  v)_{j} +1,\\
 (X\cdot w)_{i}= (X\cdot  v)_{i} -1, (\text{Lemma } \ref{le:non-vanishing})\\
  (X\cdot w)_{\ell}= (X\cdot  v)_{\ell}, \ell\not= i,j
\end{array}
$$

\nl \fbox{  Case 2}  $j<i$ and $(X\cdot v)_{h}\not=0$ for some least  $j<h<i$ then 

$$
\begin{array}{l}
(X\cdot w)_{j}= (X\cdot  v)_{j} +1,\\
 (X\cdot w)_{h}= (X\cdot  v)_{h} -1, \\
  (X\cdot w)_{\ell}= (X\cdot  v)_{\ell}, \ell\not= j, h
\end{array}
$$

\nl \fbox{  Case 3}  $j=i$ and $(X\cdot v)_{h}\not=0$ for some least  $i<h\leq n$ then
 
  $$
\begin{array}{l}
(X\cdot w)_{j}= (X\cdot  v)_{j} +1,\\
 (X\cdot w)_{h}= (X\cdot  v)_{h} -1, \\
  (X\cdot w)_{\ell}= (X\cdot  v)_{\ell}, \ell\not= j, h
\end{array}
$$

\nl \fbox{  Case 4}  $j=i$ and $(X\cdot v)_{h}=0$ for all  $h\leq n$ then
 
  $$
\begin{array}{l}
(X\cdot w)_{j}= (X\cdot  v)_{j} +1,\\
  (X\cdot  w)_{\ell}= (X\cdot v)_{\ell}, \ell\not=i
\end{array}
$$

As a result for $x_{i}>n$ we have 
$$
X(x_{i})\cdot u = X(n)\cdot u +(x_{i}-n)[X(n+1)\cdot u -X(n)\cdot u
$$
In particular if $u$ and $v$ satisfy the hypothesis of the Lemma, the two functions
$X(x_{i})\cdot u$ and $X(x_{i})\cdot v$ are equal for all $x_{i}\in \N$.

\end{proof}

The next claim is a refinement of Lemma \ref{le:action}  and shows that the equivalence $\equiv_{\text{gram}}$ can be computed.

\begin{proposition}
\label{pr:equivalence}
Two words $v$ and $v$ define the same mappings on the rows if and only if 
they have the same length and 
their restrictions coincide over  the subset $\{x_{1}, \ldots, x_{k})\in \N^{k} \mid 0\leq x_{i} \leq \max\{|u|,|v|\}+1, 1\leq i \leq k\}$.
\end{proposition}

\begin{proof} 
Observe that if they have the same restriction over the set 
of vectors of maximum coordinate  $\max\{|u|,|v|\}+1$ then by Lemma \ref{le:action} they have the same length $n$.
Let $F$ be the set of partial functions  $f:\{1, \ldots, n\}\mapsto \N$ and let $\text{supp}(f)$
 be their domain of definition. With $f\in F$ let  $X_{f}=(\alpha_{1},\ldots, \alpha_{k})$ be a formal vector where
$\alpha_{j}$ is a constant $f(j)$ when $j\in \text{supp}(f)$
and $\alpha_{j}$ is a variable otherwise. The set of variables is $\{x_{i_{1}}, \ldots, x_{i_{p}}\}$ 
where $i_{1}, \ldots, i_{p}$   is the ordered subset  $\{1, \ldots, n\}\setminus \text{supp}(f)$.
Given a word $u\in A^{*}$ we interpret $X_{f}\cdot u$ as 
a function of $\N^{p}$ into $\N^{k}$ in the natural way.
For $(\xi_{i_{1}}, \ldots, \xi_{i_{p}})\in \N^{p}$ the value of the function is denoted
 $X_{f}(\xi_{i_{1}}, \ldots, \xi_{i_{p}})\cdot u$.
We claim that for all $f\in F$, for all formal vectors  $X_{f}$  and all words $u,v$, it holds
\begin{eqnarray}
\label{eq:restriction}
\forall\ \xi_{i_{1}}, \ldots, \xi_{i_{p}}:  X_{f}(\xi_{i_{1}}, \ldots, \xi_{i_{p}}) \cdot u& =& X_{f}(\xi_{i_{1}}, \ldots, \xi_{i_{p}})\cdot v
\nonumber \\
&\Leftrightarrow&\nonumber \\
\forall\ \xi_{i_{1}}, \ldots, \xi_{i_{p}}\leq n+1:  X_{f}(\xi_{i_{1}}, \ldots, \xi_{i_{p}})\cdot u &= &X_{f}(\xi_{i_{1}}, \ldots, \xi_{i_{p}})\cdot v
\end{eqnarray}

\bigskip We show that the bottom statement  \ref{eq:restriction} implies 
the top statement.
We proceed by 
induction on $p$ by skipping the trivial case $p=0$ and the case $p=1$ which is covered
by Lemma \ref{le:one-variable}.  Define $f'\in F$ 
by the condition 
$$
f'(i)=
\left\{
\begin{array}{ll}
f(i)  &  \text{ if } i\in \text{supp}(f)\\
\xi_{i_{p}}& \text{ if } i=i_{p}\\
\text{undefined}  &  \text{otherwise}\end{array}
\right.
$$
 If $\xi_{i_{p}}\leq n+1$ then we have

$$
\begin{array}{ll}
&X_{f}(\xi_{i_{1}}, \ldots, \xi_{i_{p}})\cdot u\\
=&X_{f'}(\xi_{i_{1}}, \ldots, \xi_{i_{p-1}})\cdot u
\stackrel{\text{induction}}{=}X_{f'}(\xi_{i_{1}}, \ldots, \xi_{i_{p-1}})\cdot v\\
=&X_{f}(\xi_{i_{1}}, \ldots, \xi_{i_{p}})\cdot v
\end{array}
$$

Now define $g\in F$ by the condition 
$$
g(i)=
\left\{
\begin{array}{ll}
f(i)& \text{ if } i\in \text{supp}(f)\\
\xi_{i_{s}}  &\text{ if } 1\leq s<r\\
\text{undefined}  & \text{ otherwise} \end{array}
\right.
$$
Observe that $\text{supp}(g)=\{1, \ldots, n\}\setminus \{i_{r}\}$. Then we have
$$
\begin{array}{ll}
&X_{f}(\xi_{i_{1}}, \ldots, \xi_{i_{p}})\cdot u\\
=&X_{g}(\xi_{i_{p}})\cdot u
\stackrel{\text{Lemma \ref{le:one-variable}}}{=}X_{f'}(\xi_{i_{p}})\cdot v\\
=&X_{f}(\xi_{i_{1}}, \ldots, \xi_{i_{p}})\cdot v
\end{array}
$$

\end{proof}

The bound $n+1$ is not sharp and could probably be improved by using Lemma \ref{le:one-variable}
at the price of a more confuse statement for Proposition \ref{pr:equivalence}.

\section{Congruence in the case of three generators}

We specialize Lemma \ref{le:action} 
to three generators. The set of rows is identified with the set of triples $(x_{1},x_{2},x_{3})\in \N^{3}$
and the alphabet is renamed as $\{1,2,3\}$

\begin{equation}
\label{eq:rules}
\begin{array}{lll}
(x_{1},x_{2},x_{3})\cdot 1
&=&\left\{
\begin{array}{ll}
(x_{1}+1, x_{2}-1,x_{3}) & \text{ if } x_{2}>0\\
(x_{1}+1, x_{2},x_{3}-1) & \text{ if } x_{2}=0 \text{ and } x_{3}>0\\
(x_{1}+1, x_{2},x_{3}) & \text{ otherwise}
\end{array}
\right.\\
(x_{1},x_{2},x_{3})\cdot 2 &= &(x_{1},x_{2}+1,\max\{0,x_{3}-1\})\\
(x_{1},x_{2},x_{3})\cdot 3 &= & (x_{1},x_{2},x_{3}+1)
\end{array}
\end{equation}

\subsection{Characterization of pairs of congruent words}

Because of Corollary \ref{cor:two-letters}  we may concentrate on the pairs of words  
given by Corollary \ref{cor:a+c}.

\begin{proposition}
\label{pr:cns}
With the notations of Corollary \ref{cor:a+c}, two words $u$ and $v$ are 
$\equiv_{\text{gram}}$-congruent if and only if $c,c'\leq e$.
\end{proposition}

\begin{proof}

We compute the action   of $3^{a}2^{b}3^{c}1^{d}2^{e}3^{f}$ on $(x_{1},x_{2},x_{3})\in \N^3$
 with $a\leq b\leq d$ and $b+c\leq d+e+f$. 
Observe that 
$$
(x_{1},x_{2},x_{3})\cdot 3^{a}2^{b}3^{c}1^{d}2^{e}3^{f} = (x_{1},x_{2},x_{3})\cdot 3^{a}2^{b}3^{c}1^{d}2^{e} + (0,0,f)
$$
Thus we may suppose $f=0$ and it thus suffices to compute $(x_{1},x_{2},x_{3})\cdot 3^{a}2^{b}3^{c}1^{d}2^{e}$.
Because of Lemma \ref{cor:a+c} it all amounts to determine under which conditions, for fixed 
$b,d,e$, there exist different pairs $(a,c)$ which define the same mappings on the rows. We compute
via the rules (\ref{eq:rules})
$$
\begin{array}{rl}
& (x_{1},x_{2},x_{3}) \\
\trans{3^{a}} &(x_{1},x_{2},x_{3}+a)\\
\trans{2^{b}} %\\
&(x_{1},x_{2}+b, \max (0,x_{3}+a-b)\\
\trans{3^{c}}%\\
& (x_{1},x_{2}+b,  \max (c,x_{3}+a-b+c))
 \end{array}
$$

\nl  \fbox{Case 1} If  $x_{2}+b\geq  d$ then 

$$ \begin{array}{rl}
& (x_{1},x_{2}+b,  \max (c,x_{3}+a-b+c))\\
\trans{1^{d}} &
 (x_{1}+d,x_{2}+b-d, \max (c,x_{3}+a-b+c))\\
\trans{2^{e}}&
 (x_{1}+d,x_{2}+b-d +e, \max( 0,  \max (c,x_{3}+a-b+c))-e)\\
   & =(x_{1}+d,x_{2}+b-d +e, \max( 0, c-e, x_{3}+a-b+c-e))
   \end{array}
$$

\nl \fbox{Case 2} If $x_{2}+b<  d$ then 

$$ \begin{array}{rl}
& (x_{1},x_{2}+b,  \max (c,x_{3}+a-b+c))\\
\trans{1^{d}} & 
 (x_{1}+d,0,\max(0, \max (c,x_{3}+a-b+c))-(d-(x_{2}+b)))\\
 & =(x_{1}+d,0,\max(0, c-(d-(x_{2}+b)), x_{3}+a-b+c-(d-(x_{2}+b)))\\
 &=(x_{1}+d,0,\max(0,x_{2}+ c-d +b,  x_{2}+ c-d +x_{3}+a))\\
  \trans{2^{e}} &
 (x_{1}+d,e,\max(0,\max(0,x_{2}+ c-d +b,  x_{2}+ c-d +x_{3}+a)-e)\\
 & =(x_{1}+d,e, \max(0, x_{2}+ c-d +b-e, x_{2}+ c-d +x_{3}+a-e))
 \end{array}
$$

The same computation holds with $a'$ and $c'$ substituted for 
$a$ and  $c$.
Observe that in both  cases  1 and 2 the first two components do not depend on 
$a$ and $c$. 

Consider the third component.
In expression $\max( 0, c-e, x_{3}+a-b+c-e))$ of case 1, for $x_{3}=0$, because of $a-b\leq 0$
the condition $c-e>0$ yields
 $\max( 0, c-e, x_{3}+a-b+c-e))= c-e$. Thus if 
 $c>e$ the two words are equivalent if and only if   $c=c'$ and $a=a'$. If  
 $c-e<0$ then for  whatever values for $x_{3}$ we have
 $\max( 0, c-e, x_{3}+a-b+c-e))=\max( 0, x_{3}+a-b+c-e))$.
 Then the two words are $\gram$-equivalent if and only if $e'<c$ 
 because $\max( 0, x_{3}+a-b+c-e))= \max( 0, x_{3}+a'-b+c'-e))$. 

\end{proof}

\begin{example}
$$
\begin{array}{l}
ccc/bbb/aaabbbbb \equiv_{\text{gram}} cc/bbbc/aaabbbbb \equiv_{\text{gram}}\\
 c/bbbcc/aaabbbbb 
\equiv_{\text{gram}} bbbccc/aaabbbbb
\end{array}
$$

\end{example}

Observe that for fixed values of $b,d,e,f$  the maximum  number of occurrences of
$3$ in $w=3^{a}2^{b}3^{c}1^{d}2^{e}3^{f}$  is so that $w$ is not unique in its class is $b+e+f-1$.

\subsection{Projections on subalphabets}

Given $B\subsetneq A $ we let $\pi_{b}$ denote the projection of $A^{*}$ over $B^{*}$, i.e.,
the morphism defined by $\pi_{B}(a)=a$ if $a\in B$ and $1$ if $a\not\in B$. Routine computations show that over 
$A=\{1,2,3\}$ two $\gram$-equivalent words $u$ and $v$  so are their projections over all subalphabets 
$B\subseteq A$. 
However, the converse does not hold. Indeed, we have $23311223\not\gram23331122$ but the
row normal forms of their projections over the subalphabets $\{1,2\}$, 
$\{1,3\}$ and $\{2,3\}$ are equal, respectively $21122$, $33113$ and $332223$.

With $|A|=4$ the projections of $\gram$-equivalent words are no longer necessarily 
$\gram$-equivalent. Indeed, we have  $4213\gram 2413$ \footnote{for both words the image of $(x_{1}, x_{2}, x_{3}, x_{4})$ is
$ (x_{1}+1, x_{2}, x_{3}+1, \max\{x_{4}-1, 0\}$if  $x_{3}=0$ and
$(x_{1}+1, x_{2}, x_{3}, x_{4})$ otherwise}
but $421\not\gram 241$  because the bottom rows of the tableaux   are respectively
$1$ and $14$.

%%%%%%%%%%%%%%%%%%%%%%%%%%%%%%%%%%%%%%%%%%%%%%%%%
%%%%%%%%%%%%%%%%%%%%%%%%%%%%%%%%%%%%%%%%%%%%%%%%%
\section{The Grammic monoid over $\{1,2,3\}$}
%%%%%%%%%%%%%%%%%%%%%%%%%%%%%%%%%%%%%%%%%%%%%%%%%
%%%%%%%%%%%%%%%%%%%%%%%%%%%%%%%%%%%%%%%%%%%%%%%%%

\begin{theorem}

Let $\sim$ be the relation consisting of the pairs $(xuy, xvy)\in \{0,1,2\}^*$ where 
$(u,v)$ is one of the rules in (\ref{eq:knuth}) or the new rule 
\begin{equation}
\label{eq:relation2}
 (3212 , 2132)
\end{equation}

Then $\gram$ is the transitive closure $\stel $ of the relation  $\sim$.
Furthermore, if the words are in the form of Corollary \ref{cor:a+c}, the number of applications of the rule 
(\ref{eq:relation2}) 
to prove that they are $\gram$-equivalent is equal to  $|c-c'|$. 

\end{theorem}

\begin{proof}

We have  $3212 \gram 2132$ because for both words the image of $(x_{1}, x_{2}, x_{3})$ is
$ (x_{1}+1, x_{2}+1, \max\{x_{3}-1, 0\}$ thus
$u\stel v$ implies $u\equiv_{\text{gram}} v$. We prove the converse.
By Corollary  \ref{cor:two-letters} we may assume
that $u$ and $v$ have an occurrence of each of the three letters.
We start  with two different $\gram$-congruent row normal forms
$$
w(c)= 3^{a}2^{b}3^{c}1^{d}2^{e}3^{f}\quad \text{ and }\  w(c')=3^{a'}2^{b}3^{c'}1^{d}2^{e}3^{f}
$$
with 
\begin{equation}
\label{ex:abcdef}
0\leq c, c'\leq e,  a+c=a'+c'.
\end{equation}

 Observe that
it suffices to consider the case  $c'=c+1$ since if $c'=c+t\leq e'$ then we have 
$w(c)\sim  w(c+1)\sim \cdots \sim w(c+t)$. 
We want to show that by successive substitutions of one handside of equation \ref{eq:relation2}
for the other side 
we can rewrite $w(c)$ into $w(c')$.

\bigskip We proceed by case study where Case 1 assumes  $ d\leq b+c+1$ and 
 Case 2 assumes  $d> b+ c+1$.
We  implicitly and repeatedly use the fact that two columns commute if 
and only if one is a subset of the other.

\nl Case 1. The two words are of the following 
column normal form.
$$
\begin{array}{l}
\label{eq:1short}
(321)^{\alpha} (21)^{\beta}(31)^{\gamma}(32)^{\delta}2^{\epsilon}3^{f}\quad \text{ and }\  (321)^{\alpha-1}  (21)^{\beta+1}(31)^{\gamma}(32)^{\delta+1}2^{\epsilon-1}3^{f} \quad \epsilon>\gamma\\
\end{array}
$$

We may factorize the words respectively   as

$$
\left.
\begin{array}{l}
(321)^{\alpha} \cdot (321)(21)^{\beta}(31)^{\gamma}(32)^{\delta}2 \cdot 2^{\epsilon-1}3^{f}\\
 (321)^{\alpha} \cdot (21)^{\beta+1}(31)^{\gamma}(32)^{\delta+1}\cdot 2^{\epsilon-1}3^{f}
\end{array}
\right\}\text{ with } \epsilon \geq \gamma+1, \delta \geq 0
$$

Therefore it suffices to show how to pass from $w_{1}=(321)(21)^{\beta}(31)^{\gamma}(32)^{\delta}2^{\gamma+1}$ to $w_{2}=(21)^{\beta+1}(31)^{\gamma}(32)^{\delta+1}2^{\gamma}$.
This is obtained by applying the elementary rules concerning the commutation of two columns and 
equalities of the form $1^{\gamma} (32)^{\gamma}= (31)^{\gamma} 2^{\gamma}$

\bigskip 

$$
\begin{array}{l}
(321)(21)^{\beta}(31)^{\gamma}(32)^{\delta}2^{\gamma} 2\\
\equiv  (321)(21)^{\beta}1^{\gamma} (32)^{\delta+\gamma} 2\\
\equiv (21)^{\beta} 1^{\gamma} (321) (32)^{\delta+\gamma}2\\
\equiv (21)^{\beta} 1^{\gamma} (321) 2(32)^{\delta+\gamma}\\
\sim (21)^{\beta} 1^{\gamma} (21)(32)(32)^{\delta+\gamma}\\
\equiv (21)^{\beta+1} 1^{\gamma} (32)^{\delta+\gamma +1}\\
\equiv (21)^{\beta+1} (31)^{\gamma} 2^{\gamma} (32)^{\delta+1}\\
\equiv (21)^{\beta+1} (31)^{\gamma}  (32)^{\delta+1} 2^{\gamma}
\end{array}
$$
thus $w_{1}\stel w_{2}$. Furthermore the computation used the substitution $3212 \sim (21)(32)$
on a unique occurrence  which means that there exist two word $u,v$ such that
$$
w_{1}\equiv u3212v \sim u2312v  \equiv w_{2}
$$

\nl Case 2  The two words are of the following form.
$$
\left.
\begin{array}{l}
\label{eq:1long}
(321)^{\alpha} (21)^{\beta} (31)^{\gamma}1^{\delta}2^{\epsilon}3^{f}\\
  (321)^{\alpha-1} (21)^{\beta}(31)^{\gamma+1}1^{\delta}2^{\epsilon}3^{f} 
\end{array}
\right\} \quad \epsilon \geq \gamma +1, \delta>0
$$
and as above is suffices to show how to pass from $w_{1}=(321) (21)^{\beta} (31)^{\gamma}1^{\delta} 2^{\gamma+1}$ to $w_{2}=(21)^{\beta+1}(31)^{\gamma+1}1^{\delta-1}2^{\gamma+1}$.
This results from the following boring computation  

$$
\begin{array}{l}
(321)(21)^{\beta} (31)^{\gamma}1^{\delta} 2^{\gamma+1}\\
\equiv (321)(21)^{\beta}1^{\delta}  (31)^{\gamma} 2^{\gamma}2\\
\equiv (321)(21)^{\beta}1^{\gamma+\delta}  (32)^{\gamma}2\\
\equiv (21)^{\beta}1^{\gamma+\delta} (321) 2 (32)^{\gamma}\\
\sim (21)^{\beta}1^{\gamma+\delta} (21)(32) (32)^{\gamma}\\
\equiv (21)^{\beta+1}1^{\gamma+\delta}  (32)^{\gamma+1}\\
\equiv (21)^{\beta+1}1^{\delta}  (31)^{\gamma+1}2^{\gamma} 2\\
\equiv (21)^{\beta+1} (31)^{\gamma+1} 1^{\delta}  2^{\gamma+1}
\end{array}
$$

As in the previous case we have $w_{1}\stel w_{2}$ and  the computation used the substitution $3212 \sim 2132$
on a unique occurrence. 

\section{Further works}

Our characterization of $\gram$-congruent words over a three-letter alphabet allows one to answer some questions concerning arbitrary alphabets. 
E.g., in \cite{AR2022} the authors show the following property for  all alphabets $B\subseteq A$: if two words in $B^{*}$
have the same action of all columns labelled by $B$ they have the same action on  all columns labeled by $A$. It is routine to check that the result holds when  $|B|\leq 3$ and when ``row'' is substituted for ``column''.

A modest amount of computation led us to conjecture that in the case of four letters, the congruence is generated by the Knuth rules plus the rule $dbac=badc$ with $a<b\leq c<d$ which could probably be proved or disproved  by case study performed by a theorem prover. However, it is 
not clear how far a brute force approach could be drawn and  a more conceptual  approach is desirable.

\end{proof}

\end{document}